\documentclass[10pt]{ieeeconf}
\pdfoutput=1

\IEEEoverridecommandlockouts                              
\usepackage{graphicx}
\usepackage{graphics}
\usepackage{amsmath}
\usepackage{amssymb}
\usepackage{textcomp}
\usepackage{color}
\usepackage{caption}
\usepackage{subcaption}

\newcommand{\diag}{\operatorname{diag}}

\newcommand{\reals}{ {\mathbb{R}}}

\newcommand{\Gin} {A_{\text{ineq}}}
\newcommand{\Gcs} {c}
\newcommand{\Geq} {A_{\text{eq}}}
\newcommand{\Ceq} {b_{\text{eq}}}
\newcommand{\Cin} {b_{\text{ineq}}}

\newcommand{\Uin} {U_{\text{ineq}}}
\newcommand{\Ueq} {U_{\text{eq}}}
\newcommand{\Ucs} {U_{\text{cost}}}
\newcommand{\Ics} {I_{\text{cost}}}

\newcommand{\Ieq} {I_{\text{eq}}}
\newcommand{\Iin} {I_{\text{ineq}}}
\newcommand{\IneqSet} {\mathcal{I}}

\newcommand{\di} [1] { \diag\left({\bf 1}^T #1\right)}

\newtheorem{lemma}{Lemma}
\newtheorem{theorem}{Theorem}
\newtheorem{proposition}{Proposition}
\newtheorem{assumption}{Assumption}

\newtheorem{remark}{Remark}

\begin {document}

\title{ A Novel Method of Solving Linear Programs with an Analog Circuit}

\author{Sergey Vichik${}^{\dag}$, Francesco Borrelli${}^{\dag}$
\thanks{\textbf{${}^{\dag}$  Department of Mechanical Engineering, University of California, Berkeley, 94720-1740, USA,  \{sergv,fborrelli\}@berkeley.edu}} \\ 
31 January, 2013}
\maketitle

\begin{abstract}
We present the design of an analog circuit which solves linear programming (LP) problems.
In particular, the steady-state circuit voltages are the components of the LP optimal solution.
The paper shows how to construct the circuit and provides a proof of equivalence between the circuit
and the LP problem.
The proposed method is used to implement a LP-based Model Predictive Controller by using an analog circuit.
Simulative and experimental results show the effectiveness of the proposed approach.
\end{abstract}


\section{Introduction}
Analog circuits for solving optimization problems have been
extensively studied in the past~\cite{Dennis,Hopfield86,Chua88}.
Our renewed interests stems from Model Predictive Control (MPC)~\cite{garcia1989model},~\cite{mayne2000constrained}.
%
In MPC at each sampling
time, starting at the current state, an open-loop optimal control problem is solved over
a finite horizon. The optimal command signal is applied to the process only during the
following sampling interval. At the next time step a new optimal control problem based
on new measurements of the state is solved over a shifted horizon. The optimal solution
relies on a dynamic model of the process, respects input and output constraints,
and minimizes a performance index.
When the model is linear and the performance index is based on one-norm or $\infty$-norm, the
resulting optimization problem can be cast as a linear program (LP),
where the state enters the right hand side (rhs) of the LP constraints.

We present the design of an analog circuit whose steady state voltages
are the LP optimizers.
Thevenin Theorem is used to prove that the proposed design yields a passive circuit.
Passivity and KKT conditions of a tailored Quadratic Program are used to prove that the analog circuit 
solves the associated LP.
The proposed analog circuit can be used to repeatedly solve LPs with varying rhs and therefore
is suited for linear MPC controller implementation.
For some classes of applications the suggested implementation can be faster, cheaper and consume less power than digital implementation.
A comparison to existing literature reveals that the proposed circuit is simpler and faster than previously published designs.

The paper is organized as follows. Existing  literature is discussed in
section~\ref{sec:prev}.
We show how to construct an analog circuit from a given LP in
section~\ref{sec:circuit}. Section~\ref{sec:analys} proves the equivalence between the LP and the circuit.
Simulative and experimental results show the effectiveness of the approach in section~\ref{sec:examp}.
Concluding remarks are presented in section~\ref{sec:concl}.

\section{Previous works}
\label{sec:prev}
\subsection{Optimization problems and electrical networks}

Consider the linear programming (LP) problem
\begin{subequations}
\label{eq:genopt}
\begin{align}
    \min_{ V = [V_1,\ldots,V_n]}  & \Gcs^T  V  \label{eq:genoptcost} \\
    \text{s.t.} \hspace{0.5cm}   \Geq&  V = \Ceq   \label{eq:genopteq} \\
      \Gin&  V  \leq \Cin   \label{eq:genoptineq}
\end{align}
\end{subequations}
where $[V_1,\ldots,V_n]$ are the optimization variables, $\Gin $ and $\Geq $ are matrices, and $\Gcs$, $\Ceq $ and $\Cin $ are column vectors.

The monogram by J. Dennis~\cite{Dennis} from 1959 presents
an analog electrical network for solving the LP~\eqref{eq:genopt}.
In Dennis's work the primal and dual
optimization variables are represented by the circuit currents and voltages, respectively.
A basic version of Dennis's circuit consists of resistors, current sources, voltage
sources and diodes. In this circuit each entry of matrices $\Gin$ and $\Geq$ is equal to number of wires that are connected to a common node. Therefore, this circuit is limited to problems where the matrices $\Gin$ and $\Geq$  contain only small integer values.
An extended version of the circuit includes multiport DC-DC transformer and can represent arbitrary  matrices $\Gin$ and $\Geq$.
Current distribution laws in electrical networks
(also known as minimum dissipation of energy principle or Kirchoff's laws)
are used to prove that the circuit converges to the solution of the optimization
problem.
This work had limited practical impact due to difficulties in
implementing the circuit, and especially in implementing the multiport DC-DC transformer.

In later work, Chua~\cite{Chua1982} showed a different and more practical way to realize the multiport DC-DC
transformer using operational amplifiers. In subsequent works,
Chua~\cite{Chua88},~\cite{Chua84} and Hopefield~\cite{Hopfield86}
proposed circuits to solve non-linear optimization problem of the form
\begin{align}
\label{eq:gennlopt}
\min_x & \ f(x) \notag \\
\text{s.t.} &  \ g_j(x) \leq 0 , \ j=1 \ldots m
\end{align}
where $x \in \reals^n$ is vector of optimization variables, $f(x)$ is the cost function and $g_j(x)$ are the $m$ constraint functions.
The LP~(\ref{eq:genopt}) was solved as a special case of problem~(\ref{eq:gennlopt})~\cite{Chua88},~\cite{Hopfield86}.
The circuits proposed by Chua, Hopefield and coauthors model the Karush-Kuhn-Tucker
(KKT) conditions by representing primal variables as capacitor voltages and dual variables
as currents. The dual variables are driven by the inequality
constraint violations using  high gain amplifiers.
The circuit is constructed in a way that capacitors are charged with a current proportional to the gradient of the Lagrangian of problem~\eqref{eq:gennlopt}
\begin{equation}
\label{eq:chuaLagr}
\frac{\partial x_i} {\partial t} =  - \left [ \frac{\partial f(x)} {\partial
  x_i}  + \sum_{j=1}^m I_j \frac{\partial g_j(x)} {\partial
  x_i} \right ]
\end{equation}
where $\frac{\partial x_i} {\partial t}$ is the capacitor voltage derivative and $I_j$ is the current corresponding to the $j$-th dual variable. The derivatives $\frac{\partial f} {\partial x_i}$ and $\frac{\partial g_j} {\partial x_i}$ are implemented by using combinations of analog electrical devices~\cite{Jackson1960}.
When the circuit reaches an equilibrium, the capacitor charge is constant ($\frac{\partial x_i} {\partial t}=0$) and equation~\eqref{eq:chuaLagr} becomes one of the KKT conditions. The authors prove that their circuit always reaches an equilibrium point that satisfies the KKT conditions.
This is an elegant approach since the circuit can be intuitively mapped to the
KKT equations. However, the time required for the capacitors to reach an equilibrium is non-negligible.
This might be the reason for relatively large settling time  reported to be "tens of milliseconds" for those circuits in~\cite{Chua88}.

\subsection{Applying analog circuits to MPC problems}
The analog computing era declined before the widespread use of Model Predictive Control.
For this reason, the study of analog circuits to implement MPC problems has never been pursued.
In~\cite{Humbert2001} fast analog PI controllers are implemented on an Anadigm's Field Programmable Analog Array (FPAA) device~\cite{anadigm}
for an application involving fast chemical microreactor. The analog circuit designed
in~\cite{Humbert2001} has a computation time faster than a digital controller implementing the PI controller.
The article briefly proposes to use FPAA for MPC without specifying details.
To the best of authors knowledge, no further work has been published in this direction.

\section{Electric circuit for solving linear optimization problem}
\label{sec:circuit}
Without loss of generality, we assume that $\Gin $, $\Geq $  and  $\Gcs$   have non-negative entries.
Any LP may be transformed into this form by using a three-step procedure. First, defining a new negative and positive variable for each original variable $ V^- + V^+ = 0$, second splitting $\Gin $, $\Geq $  and  $\Gcs$ into positive and negative parts
($\Gin=\Gin^+-\Gin^- $, $\Geq=\Geq^+-\Geq^-$  and  $\Gcs=\Gcs^+-\Gcs^-$), and third
replacing $\Gin V$, $\Geq V $  and  $\Gcs^TV$ with $\Gin^+V^+-\Gin^-V^-$, $\Geq^+V^+-\Geq^-V^-$  and  ${\Gcs^+}^TV^+-{\Gcs^-}^TV^-$,
respectively.
~\\

In the beginning of this section we present the basic building blocks
which will be lately used to create a circuit that solves problem~\eqref{eq:genopt}.
The first basic block enforces equality constraints of the form~\eqref{eq:genopteq}.
The second building block enforces inequality constraints of the form~\eqref{eq:genoptineq}. The last basic block implements the cost function.

\subsection{Equality constraint}
\label{ssec:eq}
\begin{figure} [tb]
\begin{minipage}[b]{0.4\columnwidth}
\centering
\includegraphics[width=1\textwidth]{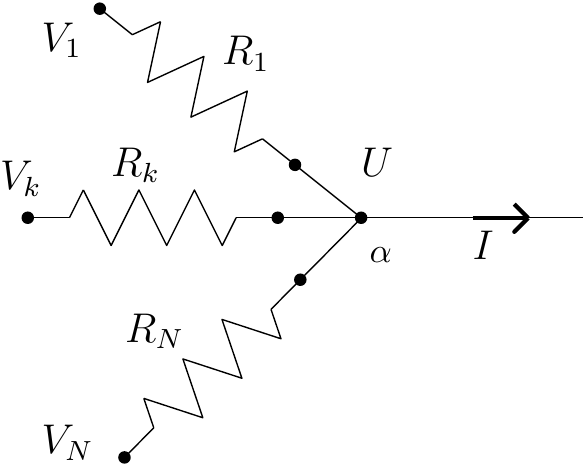}
\caption{A node with $k$ connected wires}
\label{fig:node}
\end{minipage}%
~
\begin{minipage}[b]{0.6\columnwidth}
\centering
\includegraphics[width=1\textwidth]{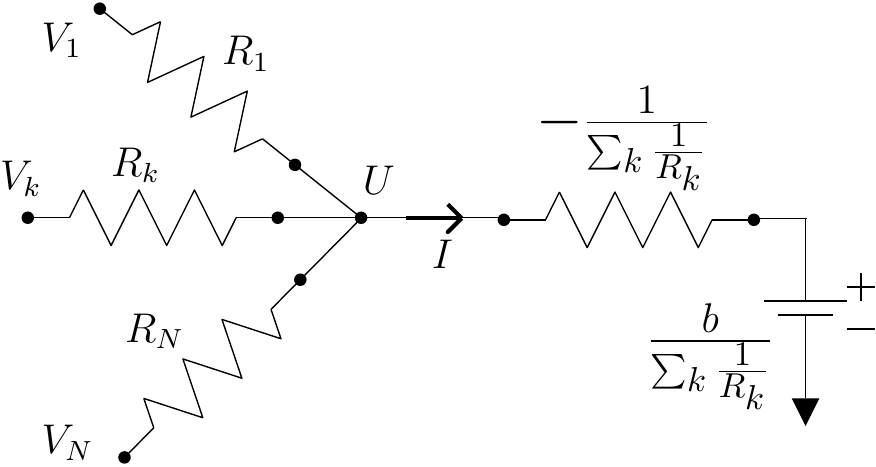}
\caption{Equality enforcing circuit. Consists of $n$ resistors $R_k$, a negative resistance and a reference voltage. }
\label{fig:eqKCL}
\end{minipage}%
\end{figure}

Consider the circuit depicted in Fig.~\ref{fig:node}.
In this circuit $n$ wires are
connected to a common node. We call this node $\alpha$, its  potential is $U$ and the current that exits this node is $I$.  Kirchhoff's current law (KCL) implies
\begin{equation}
\label{eq:KCL}
\sum_{k=1}^n I_k = \sum_{k=1}^n \frac{V_k-U}{R_k} = I,
\end{equation}
where $V_k$ is the potential of node $k$, $R_k$ is the resistance between node $k$ and the node $\alpha$.
Equation~\eqref{eq:KCL} can be written as an equality constraint on potentials $V_k$:
\begin{equation}
\label{eq:Vksum1}
\sum_{k=1}^n \frac{V_k}{R_k} = I + U  \sum_{k=1}^n \frac{1}{R_k}.
\end{equation}
If we can set the right hand side (rhs) of~\eqref{eq:Vksum1} to any desired value $b$, then~\eqref{eq:Vksum1} enforces an equality constraint on a linear combinations of $V_k$. Therefore every equality constraint~\eqref{eq:genopteq} can be implemented with a circuit which enforces~\eqref{eq:Vksum1} and  implements
\begin{equation}
\label{eq:V0law}
U  = \frac{b-I}{\sum_{k=1}^n \frac{1}{R_k}}.
\end{equation}
Equation~\eqref{eq:V0law} together with~\eqref{eq:Vksum1}  yields
\begin{equation}
\label{eq:Vkeq}
 \left [ \frac{1}{R_1}   \   \ldots \  \frac{1}{R_n} \right ] \left [ \begin{array}{c} V_1 \\ \vdots \\ V_n \end{array} \right ] = b .
\end{equation}
and the circuit implementing~\eqref{eq:Vkeq} is shown in Fig.~\ref{fig:eqKCL}.

\begin{remark}
In the circuit in Fig.~\ref{fig:eqKCL} the negative resistance
$-\frac{1}{\sum_{k} \frac{1}{R_k}}$ can be realized by using operational amplifiers.
\end{remark}

 \subsection{Inequality  constraint}
\label{ssec:ineq}
\begin{figure} [tb]
\begin{minipage}[b]{0.55\columnwidth}
\centering
\includegraphics[width=1\textwidth]{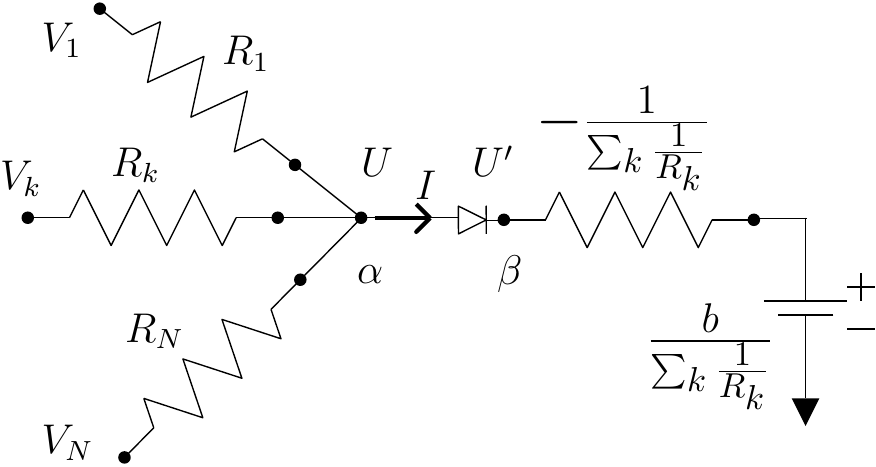}
\caption{Inequality enforcing circuit. }
\label{fig:ineqnode}
\end{minipage}
\begin{minipage}[b]{0.44\columnwidth}
\centering
\includegraphics[width=1\textwidth]{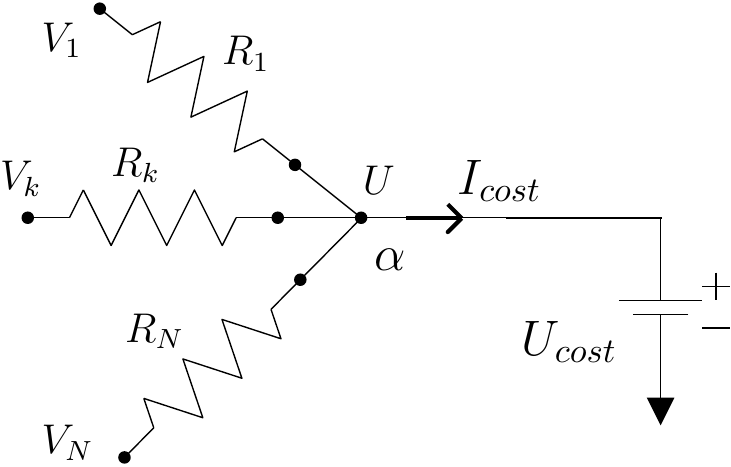}
\caption{Cost circuit }
\label{fig:costnode}
\end{minipage}
\end{figure}

Consider the circuit shown in Fig.~\ref{fig:ineqnode}.
Similarly to the equality constraint circuit, $n$ wires are
connected to a common node $\alpha$. Its  potential is $U$ and the current exiting this node is $I$.
Kirchhoff's current law (KCL) implies
\begin{equation}
\label{eq:KCL2}
\sum_{k=1}^n I_k = \sum_{k=1}^n \frac{V_k-U}{R_k} = I.
\end{equation}
An ideal diode connects node $\alpha$ to node $\beta$. The potential of node $\beta$ is $U'$. The diode enforces $U \le U'$.
In Fig.~\ref{fig:ineqnode}, the voltage $U'$ can be computed as follows
\begin{equation}
\label{eq:U'law}
U'  = \frac{b-I}{\sum_{k=1}^n \frac{1}{R_k}} \geq U.
\end{equation}
Equation~\eqref{eq:KCL2} and $U \le U'$ yield
\begin{align}
\sum_{k=1}^n \frac{V_k}{R_k}  &= I + U  \sum_{k=1}^n \frac{1}{R_k} \le I + U' \sum_{k=1}^n \frac{1}{R_k} = b.
\end{align}
Which can be compactly rewritten as
\begin{equation}
\label{eq:Vkineq}
\left [ \frac{1}{R_1}   \   \ldots \  \frac{1}{R_n} \right ] \left [ \begin{array}{c} V_1 \\ \vdots \\ V_n \end{array} \right ] \leq b,
\end{equation}
with the diode enforcing
\begin{subequations}
\begin{align}
&I \ge 0 , \label{eq:diodeIpos} \\
&I (U-U') = 0. \label{eq:diodecompl0}
\end{align}
\end{subequations}
By using~\eqref{eq:U'law} and rearranging some terms, equation~\eqref{eq:diodecompl0} can be rewritten as:
\begin{align}
  I \left (  \left (\sum_{k=1}^n \frac{1}{R_k}\right )U  - b  + I \right )
= 0 .
\label{eq:diodecompl}
\end{align}

\subsection{Cost function}
\label{ssec:cost}
Consider the circuit in Fig.~\ref{fig:costnode}. In this circuit the potential of node $\alpha$ is equal to $\Ucs$ and the current that exits the node is $I_\text{cost}$.
From~\eqref{eq:Vksum1} we have
\begin{align}
\sum_{k=1}^n \frac{V_k}{R_k}  &= I_\text{cost} + \Ucs  \sum_{k=1}^n \frac{1}{R_k} \triangleq J.
\end{align}
where $c=[1/R_1\ \ldots\ 1/R_n]$ and  $J$ is the cost function.

This part of the circuit implements the minimization of the cost function. When $\Ucs$ is set to a low value, the voltages $V_k$ are driven to a direction which leads the objective function value $J$ to approach the $\Ucs$ value. However, the cost $J$ is different from $\Ucs$ because the current $\Ics$ is not zero.
A detailed explanation on this part of the circuit will be presented later in section~\ref{ssec:equiv}.

\subsection{Connecting the basic circuits}

This section presents how to construct the circuit that solves a general LP.
We construct the conductance matrix $G \in \reals ^{(m+1) \times n}$ as
\begin{align}
G \triangleq \left [ \begin{array}{c}  \Gcs^T \\  A \end{array} \right ] = \left [ \begin{array}{c}  \Gcs^T \\  \Geq \\ \Gin \end{array} \right ]
\end{align}
and denote $G_{ij}$  the $i,j$ element of $G$.
For a given LP~\eqref{eq:genopt} the $R_{ij}$ resistor is defined as
\begin{align}
R_{ij} = \frac{1}{G_{ij}},~i=0,\ldots m, j=1,\ldots,n
\label{eq:Rij}
\end{align}
where the first row of G (corresponding to $\Gcs^T$) is indexed by 0.

\begin{figure}[tb]
\centering
   \includegraphics[width=0.42\textwidth]{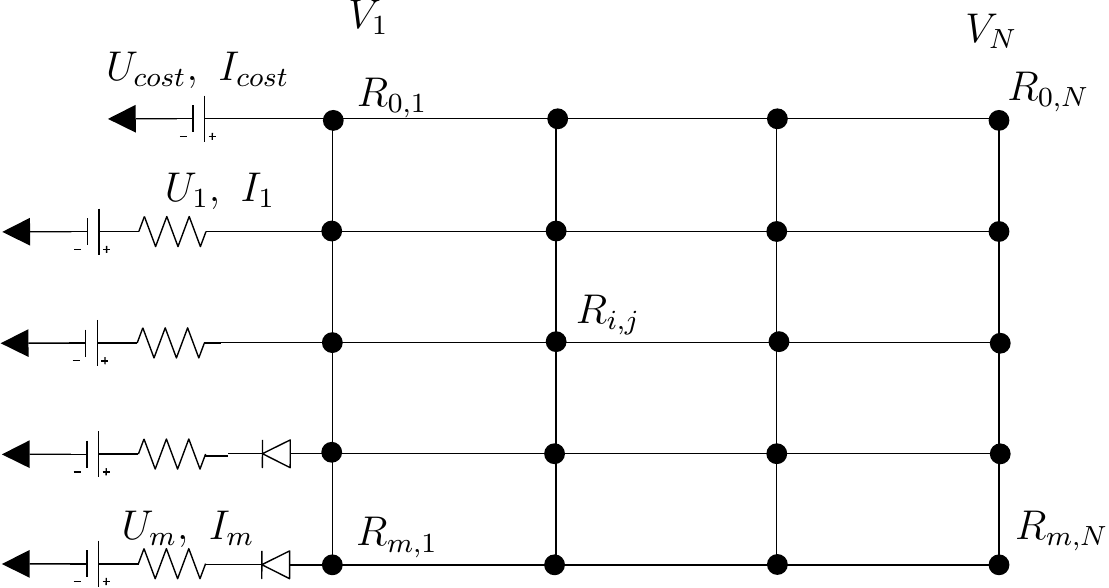}\\
  \caption{Electric Circuit solving a LP. Vertical wires are variable nodes with potentials $V_1 \ldots V_n$.  Black dots represent resistances that connects vertical and horizontal wires. Horizontal wires are cost or constraint nodes.  Each horizontal wire is connected to a ground via a negative resistance, a constant voltage source and a diode for inequalities nodes. The topmost horizontal wire is the cost circuit and is connected to a constant voltage source.}\label{fig:circgen}
\end{figure}
Consider the circuit shown in Fig.~\ref{fig:circgen}.
The circuit is shown using a compact notation where each resistor $R_{ij}$ is represented by a dot, vertical wires represent variables nodes with potentials  $V_1 \ldots V_n$ and horizontal wires represent \emph{constraint nodes}. If $G_{ij}=0$ then no resistor is present in the corresponding dot.
This circuit is constructed by connecting the  nodes associated with the variables  $V_1 \ldots V_n$ to all three types of the basic circuits: equality, inequality and cost. We will refer to such nodes as \emph{variable nodes}. Each row of the circuit in Fig.~\ref{fig:circgen} is one of the basic circuits presented in Sections~\ref{ssec:eq}, \ref{ssec:ineq} and~\ref{ssec:cost}. We claim that, if $\Ucs$ is ``small enough'', then the values of the potentials $V_1 \ldots V_n$ in this circuit are a solution of~\eqref{eq:genopt}. This claim is proven in the next section.

\begin{remark}
Some of the potentials $V_i$ may be forced externally to a desired value. By doing so, the circuit can solve different optimization problems for varying values of those potentials. This  is equivalent to adding equality constraints $V_i=b_i$ to~\eqref{eq:genopt} and modifying the value of the equality constraint free parameter $b_i$.
\end{remark}

\begin{remark}
The circuit as shown in Fig.~\ref{fig:circgen} contains no  dynamic elements such as capacitor or inductance.
Therefore, the  time required to reach steady-state is governed by the parasitic effects (e.g. wires inductance and capacitance) and by the properties of the elements used to realize negative resistance (usually opamp) and diode. Hence, a good electronic design can achieves solution times in the order of these parasitic effects. This could lead to time constants as low as a few nanoseconds.
\end{remark}

\section{Analysis of the electric circuit properties}
\label{sec:analys}

In this section we show that the circuit in Fig.~\ref{fig:circgen} with $R_{ij}$ as defined by~\eqref{eq:Rij}, is
a solution of the optimization problem~\eqref{eq:genopt} for a range of $\Ucs$ values.
First we derive the steady state equation of the electric circuit and then we show the equivalence.

\subsection{Steady state solution}

Consider the circuit in Fig.~\ref{fig:circgen}.  Let $U=[U_1,\ldots, U_m]^T$ be the voltages of the constraint nodes as shown on Fig.~\ref{fig:circgen}.  By
applying the KCL (Kirchhoff's current law) to  every variable node with potential $V_1,\ldots,V_n$ we obtain
\begin{align}
& G_{0,j}(\Ucs-V_j)+ \sum_{i=1}^{m} G_{i,j}(U_i-V_j) = 0, \hspace{0.2cm} j=1,\ldots,n
\end{align}
which can be rewritten in the matrix form
\begin{align}
 \left [ \begin{array}{ccc}
\Gcs_{1} & . & \Gcs_{n} \\
A_{11} & . & A_{1N} \\
\vdots & \cdot  & \vdots \\
A_{m1} & . & A_{mN} \end{array} \right ] ^T
\left [ \begin{array}{c}
\Ucs \\ U_1 \\ \vdots \\ U_m
 \end{array} \right ]   =
 \left [ \begin{array}{c}
 (\sum_{i=0}^{m} G_{i,1})V_1  \\  \vdots \\ (\sum_{i=0}^{m} G_{i,n})V_n
 \end{array} \right ].
\label{eq:V_KCL} \end{align}
Equation~\eqref{eq:V_KCL} can be compactly rewritten as
  \begin{align}
 & \ \Gcs\Ucs+ A^T U =
\diag(\Gcs^T+ {\bf 1}^T A)  V
\label{eq:circUV}
 \end{align}
where ${\bf 1}$ is vector of ones and $\diag(x)$ is a diagonal matrix with $x$ on its diagonal.

Next, we apply KCL on all nodes with potentials $[\Ucs, U_1,\ldots, U_m]$ to obtain
\begin{align}
& \sum_{j=1}^{n} \Gcs_{j}(\Ucs-V_j) = \Ics\\
& \sum_{j=1}^{n} G_{i,j}(U_i-V_j) = I_{i}, \hspace{0.5cm} i=1,\ldots, m
\end{align}
which  can be written in matrix form
\begin{align} \hspace{-0.2cm}
 \left [ \begin{array}{ccc}
\Gcs_{1} & . & \Gcs_{n} \\
A_{11} & . & A_{1N} \\
\vdots & \cdot  & \vdots \\
A_{m1} & . & A_{mN} \end{array} \right ]
&\left [ \begin{array}{c}
V_{1} \\ \vdots \\ V_{n}
 \end{array} \right ] 
  =  \notag \\   &
 \left [ \begin{array}{c}
 \Ucs  \sum_{j=1}^{n} \Gcs_{j} \\ U_1 \sum_{j=1}^{n} A_{1,j} \\  \vdots \\ U_m \sum_{j=1}^{n} A_{m,j}
 \end{array} \right ]   +  \left [ \begin{array}{c} \Ics \\ I \end{array} \right ].
\label{eq:U_KCL} \end{align}
Equation~\eqref{eq:U_KCL} can be compactly rewritten as
\begin{align}
  & c^TV =  {\bf 1}^T\Gcs  \Ucs + \Ics \label{eq:GVUIcost} \\
  &A V = \di{A^T}  U + I.
  \label{eq:GVUI}
\end{align}

The equality voltage regulator law~\eqref{eq:V0law} and the inequality law~\eqref{eq:U'law} can be compactly written as
\begin{subequations}
\label{eq:U_Law}
\begin{align}
&\di{\Geq^T}\Ueq = \Ceq - \Ieq \\
&\di{\Gin^T}\Uin \le \Cin - \Iin.
\end{align}
\end{subequations}
By substituting \eqref{eq:U_Law} into \eqref{eq:GVUI} we obtain
\begin{align}
& \Geq V = \Ceq \label{eq:circEQ}\\
&  \Gin V \le \Cin .\label{eq:circINEQ}
\end{align}
Substitution of~\eqref{eq:GVUI}  for inequalities to the diode constraint~\eqref{eq:diodecompl} yields
\begin{align}
\left [\Gin V - \Cin \right ]_i [\Iin]_i = 0, \ \forall i\in \IneqSet \label{eq:circDIODE}
\end{align}
where $\IneqSet$ is the set of all inequalities constraints.

We collect \eqref{eq:circUV},  \eqref{eq:GVUIcost}, \eqref{eq:GVUI}, \eqref{eq:circEQ}, \eqref{eq:circINEQ} and~\eqref{eq:diodeIpos} into one  set of equations
which characterize the circuit
\begin{subequations}
\label{eq:eleceq}
\begin{align}
& A V = \di{A^T} U + I \label{eq:fcUI} \\
&\Gcs \Ucs +  A^T U =  \diag(\Gcs^T+ {\bf 1}^T A)  V \label{eq:fcLAG} \\
& \Geq V = \Ceq \label{eq:fcEQ}\\
& \Gin V \le \Cin \label{eq:fcINEQ}\\
&\Iin \geq 0 \label{eq:fcIneg} \\
& \left [\Gin V - \Cin \right ]_i [\Iin]_i = 0, \forall i\in \IneqSet \label{eq:fcDIODE}\\
&c^TV =  {\bf 1}^T \Gcs  \Ucs + \Ics  \label{eq:fcUIcost} ,
\end{align}
\end{subequations}
where $U$, $I$, $\Ics$ and $V$ are the unknowns. The voltage of the cost node, $\Ucs$, is set externally.


\subsection{Circuit passivity}

\begin{figure} [tb]
\centering
\includegraphics[width=0.35\textwidth]{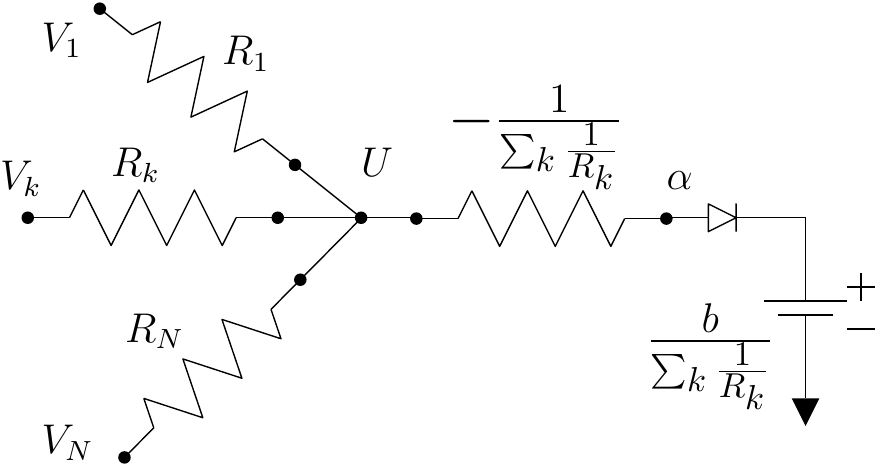}
\caption{An alternative inequality enforcing node. When diode is not present, the
circuit is an equality enforcing node.}
\label{fig:ineqnodeswapped}
\end{figure}

\begin{figure} [tb]
\centering
\includegraphics[width=0.32\textwidth]{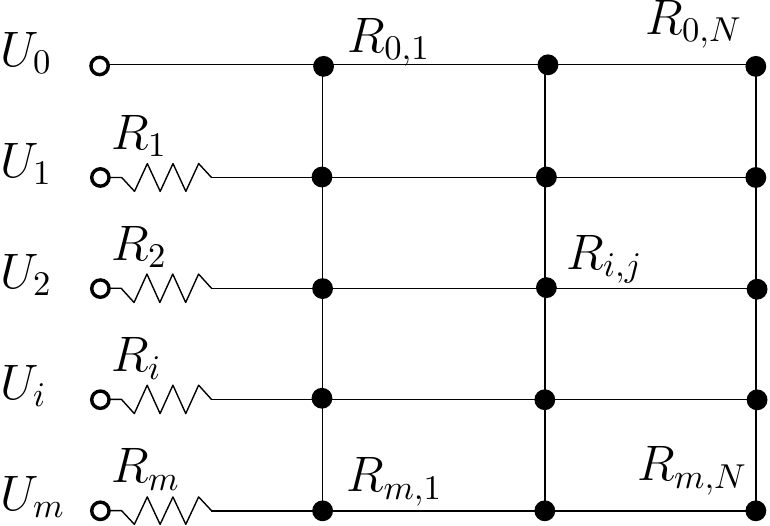}
\caption{N-port resistor network with ports $U_i$. All $R_{i,j}$ are positive resistances, all $R_k$ are negative resistances. }
\label{fig:elecnetw}
\end{figure}

We are interested in showing that the general circuit in Fig.~\ref{fig:circgen} is a passive system. 
First we swap the diode and the resistor in Fig.~\ref{fig:ineqnode} to obtain a basic inequality circuit as shown in
Fig.~\ref{fig:ineqnodeswapped}. It is simple to prove that this swap yields an equivalent
electric circuit. Next, we examine an N-port  resistor network that includes all
the resistors of the original circuit shown in Fig.~\ref{fig:circgen}, including the negative resistances, but not including the
diodes nor the constant voltage sources as shown in
Fig.~\ref{fig:elecnetw}. The ports of the electric network is
the set of all nodes marked $\alpha$ in Fig.~\ref{fig:ineqnodeswapped}.


\begin{proposition}[Network non-negativity]
\label{claim:nonneg}
The resistance network in Fig.~\ref{fig:elecnetw} is equivalent to a resistance network with non-negative
resistors.
\end{proposition}

\begin{proof}[Proof of non-negativity proposition]
Our goal is to obtain a lower bound of an equivalent resistance between any
two ports. From Fig.~\ref{fig:elecnetw} we see that a sub-network that
connects two ports consists of two negative resistances ---  one for each
port, and a mesh of positive resistors between them. We want to find an
equivalent resistance, that exist according to the Thevenin
theorem~\cite{chen2004electrical}. Let $U_i$ and $U_j$ be the two nodes in
question. Next, motivated by a fact that replacement of any of positive
resistances with a zero resistance may only reduce the total equivalent
resistance, we make a conservative assumption that all the resistors in this
network, excluding resistors directly connected to negative resistors of the
$U_i$ and $U_j$ nodes, are zero, thus $R_{k,l}=0, \forall k,l$  s.t.  $k\neq
i,j$.  In this case all variables nodes have the same
potential. This sub-network is illustrated in
Fig.~\ref{fig:reducedNetw}. The equivalent resistance of this network is
zero, since according to~\eqref{eq:V0law} the negative resistance is constructed to be equal to
the negative of parallel combination of other node resistances. For the 
$\Ucs$ cost port (which does not have the negative resistor), the
equivalent resistance is strictly greater than zero. Therefore, the equivalent
resistance between any two ports is at least zero.
\begin{figure} [tb]
\centering
\includegraphics[width=0.37\textwidth]{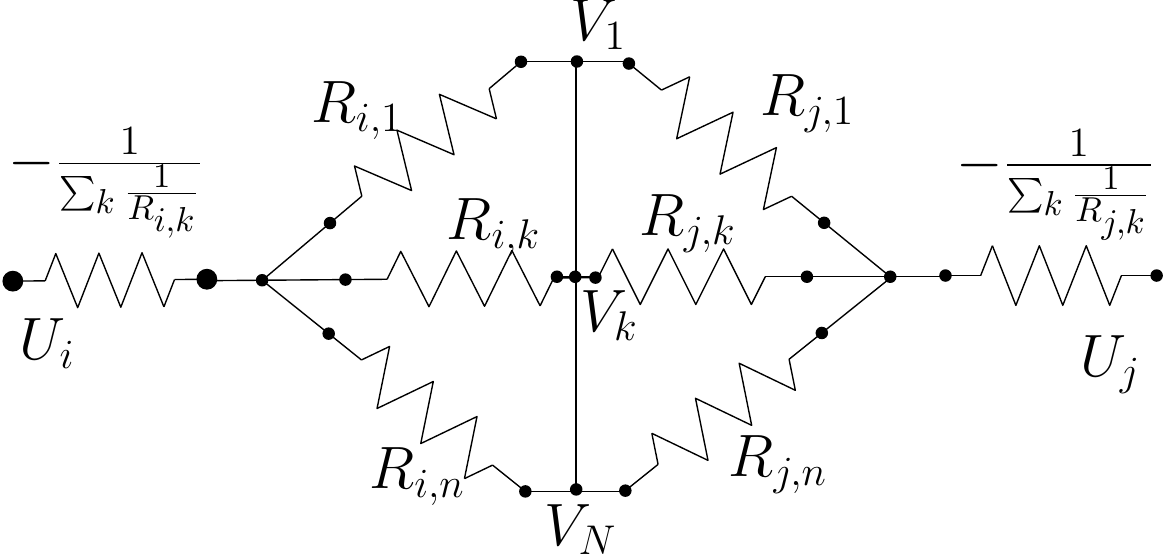}
\caption{Subnetwork that connects nodes $i$ and $j$, after assuming that all other resistors are zero.}
\label{fig:reducedNetw}
\end{figure}
\end{proof}

\subsection{Equivalence of the optimization problem and the electric circuit}
\label{ssec:equiv}
We consider the following assumptions.
\begin{assumption}
\label{ass1}
The LP~\eqref{eq:genopt} is feasible and that the set of primal optimal solutions is bounded.\end{assumption}

\begin{assumption}
\label{ass2}
The dual of LP~\eqref{eq:genopt} is feasible and that the set of dual optimal solutions is bounded.\end{assumption}

\begin{assumption}
\label{ass3}
In the LP~\eqref{eq:genopt},  $G$ is non-negative, ${\bf 1}^T G > 0 $ and ${\bf 1}^T G^T > 0 $.
\end{assumption}
\vspace{0.2cm}
\begin{theorem}
[circuit equivalence]
\label{thrm:equiv}
Let Assumptions~\ref{ass1}-\ref{ass3} hold.
Then, there exists $\Ucs^\text{crit}$, such that a solution $V^*$ to~\eqref{eq:eleceq} is also an optimizer of the LP~\eqref{eq:genopt} for all $\Ucs \leq \Ucs^\text{crit}$.
\end{theorem}

\begin{remark}
As explained earlier, the assumption on the  non-negativity of $G$ in Theorem~\ref{thrm:equiv} is not restrictive.
Also, ${\bf 1}^T G > 0 $ and ${\bf 1}^T G^T > 0 $ are always satisfied for LP problems without zero rows or zero columns.
\end{remark}
\begin{remark}
In Theorem~\ref{thrm:equiv} we require that the sets of primal optimal and dual optimal solutions are bounded. This can be guaranteed if the primal feasible set is bounded and linear independent constraint qualification (LICQ) holds.
\end{remark}

The theorem will be proven in the following way: first we claim that the equations~\eqref{eq:fcUI}-\eqref{eq:fcDIODE}  have a solution when no cost function is present ($\Gcs=0$); second, we show that there exists $\Ucs^\text{crit}$ such that any solution to~\eqref{eq:eleceq} is also an LP solution; third, we show that for all $\Ucs \leq \Ucs^\text{crit}$ any solution to~\eqref{eq:eleceq} is also an LP solution.

Consider an electric circuit, that consists of constraint sub~circuits and no cost sub~circuit.
This electric circuit is characterized by~(\ref{eq:fcUI})-(\ref{eq:fcDIODE}) with $\Gcs=0$.
\begin{lemma} [Existence of solution to a no-cost circuit]
\label{lem:existSol}
Let Assumption~\ref{ass1} hold. Assume that $A$~is non-negative, ${\bf 1}^T A > 0 $ and ${\bf 1}^T A^T > 0 $.
Then, the equations~(\ref{eq:fcUI})-(\ref{eq:fcDIODE}) have a solution when $\Gcs=0$.
\end{lemma}

\begin{proof}
First we rearrange~(\ref{eq:fcUI})-(\ref{eq:fcDIODE}).
Equation~\eqref{eq:fcUI} can be split into an equality and inequality parts
\begin{align}
& \Geq  = \di{\Geq^T} \Ueq + \Ieq \label{eq:UIeq} \\
& \Gin  = \di{\Gin^T} \Uin + \Iin \label{eq:UIin} .
\end{align}
Equation~\eqref{eq:fcLAG} can be rewritten as
\begin{align}
& \Geq^T \Ueq + \Gin^T \Uin   =  \di{A }  V \label{eq:LAG2}.
 \end{align}
Therefore,~(\ref{eq:fcUI})-(\ref{eq:fcDIODE})  can be written as
 \begin{subequations}
\label{eq:eleceqnocost}
\begin{align}
& \Geq V = \di{\Geq^T} \Ueq + \Ieq \label{eq:UIeq} \\
& \Gin V = \di{\Gin^T} \Uin + \Iin \label{eq:UIin}   \\
& \Geq^T \Ueq + \Gin^T \Uin   =  \di{A }  V \label{eq:LAG2} \\
& \Geq V = \Ceq \label{eq:fc2EQ}\\
& \Gin V \le \Cin \label{eq:fc2INEQ}\\
& \Iin \geq 0 \\
& \left (\Gin V - \Cin \right )_i{\Iin}_i = 0,\ \forall i\in \IneqSet \label{eq:fc2DIODE}.
\end{align}
\end{subequations}

Next, consider the following quadratic program (QP)
\begin{subequations}
\label{eq:primfeas}
 \begin{align}
 \min_V  &\ V^T Q V \notag \\
s.t. \ &  \Geq V=\Ceq  \\
 &\Gin V \le \Cin,
 \end{align}
\end{subequations}
This problem has a finite solution for any $Q$ because the feasibility domain is bounded and not empty. The value of $Q$ will be selected later.
We use this problem to find a solution to~(\ref{eq:fcUI})-(\ref{eq:fcDIODE}). KKT is a necessary optimality condition for problems with linear constraints (Theorem~5.1.3 in~\cite{bazaraa2006nonlinear}), therefore, there exist~$V^\star$, $\mu^\star$, $\lambda^\star$ which satisfy the KKT conditions
\begin{subequations}
\label{eq:KKTfprob}
\begin{align}
&\Geq^T \mu^\star + \Gin^T   \lambda^\star  + QV^\star = 0 \\
     &  \Geq   V^\star = \Ceq \\
     &  \Gin V^\star  \leq \Cin  \\
& \lambda^\star \geq 0 \\
&   (\Gin V^\star  - \Cin)_i \lambda_i^\star=0, \ i \in \IneqSet,
\end{align}
\end{subequations}
where $\mu^\star$ and $\lambda^\star$ are the dual variables.

We choose $Q$ and use $\mu^\star$, $\lambda^\star$ and $V^\star$ to compute $\Ueq^\star$, $\Uin^\star$, $\Ieq^\star$ and $\Iin^\star$
 \begin{subequations}
 \label{eq:QIUdef}
 \begin{align}
 Q= &\di{A }   -\Geq^T\di{\Geq^T}^{-1}\Geq  \notag \\ & -\Gin^T\di{\Gin^T}^{-1}\Gin \\
\Ieq^\star =&   \di{\Geq^T}\mu^\star \\
\Ueq^\star = &   \di{\Geq^T}^{-1}\Geq V^\star - \mu^\star  \\
\Iin^\star = &\di{\Gin^T}\lambda^\star   \\
\Uin^\star = &  \di{\Gin^T}^{-1}\Gin V^\star  - \lambda^\star \label{eq:lambdadef}.
 \end{align}
 \end{subequations}
Note that $\di{\Gin^T}$ and $\di{\Gin^T}$ are invertible and positive from the assumptions of Lemma~\ref{lem:existSol}. Equations~\eqref{eq:QIUdef} are combined with~\eqref{eq:KKTfprob} to get
\begin{subequations}
\label{eq:feasprimdual2}
 \begin{align}
 & \Geq V^\star = \di{\Geq^T} \Ueq^\star + \Ieq^\star \label{eq:UIeq3}  \\
  & \Gin V^\star = \di{\Gin^T} \Uin^\star + \Iin^\star \label{eq:UIin3}  \\
& \Geq^T \Ueq^\star + \Gin^T \Uin^\star  =  \di{A }  V^\star \label{eq:feaspd1}   \\
&  \Geq V^\star=\Ceq  \\
 &\Gin V^\star \le \Cin \\
 & \Iin^\star \geq 0 \\
&   (\Gin V^\star  - \Cin)_i {\Iin}_i^\star=0, \ i \in \IneqSet.
\end{align}
\end{subequations}
Equations~\eqref{eq:feasprimdual2} have a solution and are identical to~\eqref{eq:eleceqnocost}.
Therefore, there exist $V^\star$, $U^\star$ and $I^\star$ solving~(\ref{eq:fcUI})-(\ref{eq:fcDIODE}) when $\Gcs=0$.
\end{proof}

Our next goal is to show that there exists a $\Ucs$ such that circuit solution is also an LP~\eqref{eq:genopt} solution. To show this we concatenate the primal problem \eqref{eq:genopt} with a corresponding dual problem~\cite{bertsimas1997introduction}
\begin{subequations}
\label{eq:gendual}
\begin{align}
\max_\lambda & \ b^T \lambda \\
\text{s.t.} &\ [\Geq^T\ \Gin^T] \lambda = \Gcs \\
&\ \left [0 \  I_{|\IneqSet|}  \right ] \lambda \geq 0,
\end{align}
\end{subequations}
where $I_{|\IneqSet|}$ is an identity matrix of size equals to number of inequality constraints.
We create the following feasibility problem
\begin{subequations}
\label{eq:primdual}
\begin{align}
\min_{\lambda,V} & \ 0 \\
    \text{s.t.}& \ \Geq   V = \Ceq, \
       \Gin V  \leq \Cin  \\
&\ [\Geq^T\ \Gin^T] \lambda = \Gcs, \
 \left [0 \   I_{|\IneqSet|}  \right ] \lambda \geq 0\\
&\    \Gcs^T V + b_-^T \lambda + b_+^T \lambda_-  = 0, \ \lambda +  \lambda_- =0,\label{eq:pdslackness}
\end{align}
\end{subequations}

where $b_+$ and $b_-$ are the absolute values of the positive and the negative components of $b$ and $\lambda_-$ equals to $-\lambda$. Note that \eqref{eq:pdslackness} is equivalent to $\Gcs^T V =  b^T \lambda$.
\begin{remark}
\label{rem:Apos}
From the Assumption~\ref{ass3} and from the structure of~\eqref{eq:pdslackness}, it follows that the matrix of equality and inequality constraints has non-negative coefficients and non-zero rows and columns.
\end{remark}
All feasible points of problem~\eqref{eq:primdual} are primal~\eqref{eq:genopt} and dual~\eqref{eq:gendual} optimal
solutions~\cite{bertsimas1997introduction}.

\begin{figure} [tb]
\centering
\includegraphics[width=0.42\textwidth]{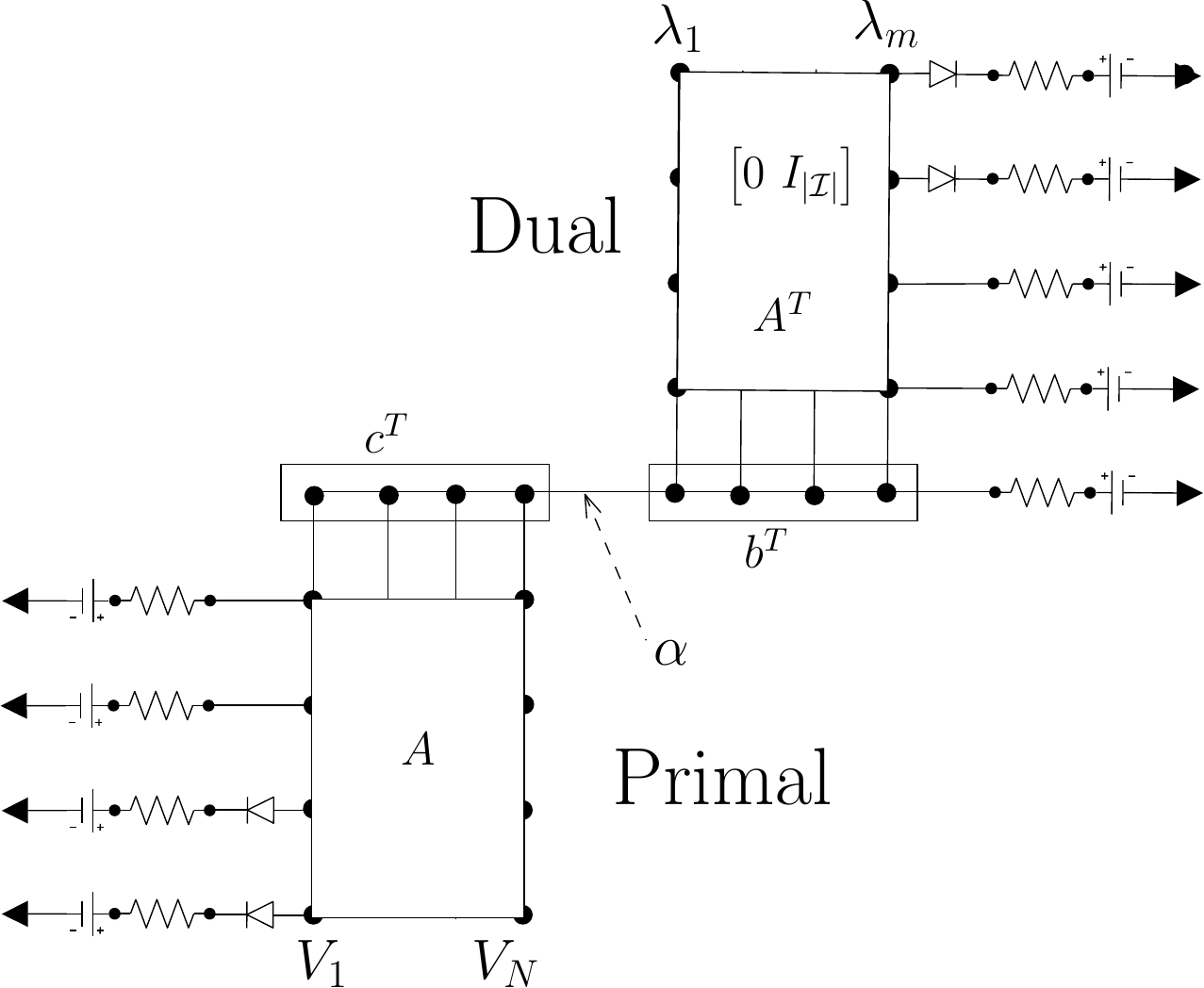}
\caption{Circuit that implements the primal-dual feasibility problem. Primal and dual constraints are separated. Primal and dual parts are connected only via zero duality gap constraint.  For compactness, $b_+$ and $b_-$ are represented as $b$ and $\lambda_-$ is embedded in $\lambda$. }
\label{fig:primdual}
\end{figure}

Problem~\eqref{eq:primdual} is solved by the circuit shown in Fig.~\ref{fig:primdual}. This circuit contains two parts: the primal and the dual circuits, each has the general form as in Fig.~\ref{fig:circgen} and consists of equality and inequality sub~circuits corresponding to constraints of the primal and dual problems. Note that no cost circuit is present in the primal and in the dual circuit. Instead, those circuits are connected by equality sub~circuit that corresponds to the zero duality gap constraint~\eqref{eq:pdslackness}.
\begin{proposition}
Let Assumptions~\ref{ass1}-\ref{ass3} hold.
The circuit in Fig.~\ref{fig:primdual} admits a solution. Moreover, for any circuit solution, the voltages $V$ of the variable nodes are a solution to the original LP~\eqref{eq:genopt}.
\end{proposition}
\begin{proof}
The circuit in Fig.~\ref{fig:primdual} consists only of equality and inequality sub~circuits. As shown in sections~\ref{ssec:eq} and~\ref{ssec:ineq} the   variable nodes voltages must satisfy the associated equality or inequality constraints and thus equations~\eqref{eq:primdual}.  The feasible set of problem~\eqref{eq:primdual} is the set of all primal optimal and dual optimal variables of   problem~\eqref{eq:genopt}.
This feasible set is bounded by assumption.
This fact and the results from Remark~\ref{rem:Apos} imply that all the assumptions of Lemma~\ref{lem:existSol} are satisfied. We conclude that the circuit admits a solution. Moreover, every solution must be a solution of the original LP~\eqref{eq:genopt}, because it satisfies simultaneously dual and primal problems with zero duality gap~\cite{bertsimas1997introduction}.
\end{proof}

In the circuit shown in Fig.~\ref{fig:primdual}, the dual and the primal circuits are connected with a single wire that has some voltage when the circuit settles. We call this voltage $\Ucs^{crit}$.
\begin{lemma}[Exists $\Ucs^{crit}$ ]
\label{lem:existsUc}
Let Assumptions~\ref{ass1}-\ref{ass3} hold.
Consider the circuit in Fig.~\ref{fig:circgen} and its corresponding equations~\eqref{eq:eleceq}.
  A solution $V^*$ to~\eqref{eq:eleceq} with
$\Ucs = \Ucs^\text{crit}$ is an optimizer of the LP~\eqref{eq:genopt}.
\end{lemma}
\begin{proof}
 If  a voltage equals to $\Ucs^\text{crit}$ is applied externally to the wire that connects the primal and the dual parts (at point $\alpha$ in Fig.~\ref{fig:primdual}), we can remove the dual circuit without affecting the primal one.
Therefore, the circuit in Fig.~\ref{fig:circgen} admits the same solution as the primal circuit in Fig.~\ref{fig:primdual}.
\end{proof}


To complete the proof of Theorem~\ref{thrm:equiv} we need to show that for any voltage $\Ucs \leq \Ucs^{crit}$ the circuit will continue to yield the optimal solution.
Assume that  $\Ucs$ is perturbed by $ \Delta \Ucs$ from the value $\Ucs^{crit}$. We denote
perturbed values in variable voltages $V$ and the cost current $\Ics$ as $\Delta V$ and $\Delta \Ics $.
Next, we examine the Thevenin equivalent resistance~\cite{chen2004electrical}  as seen from the cost node. From  Proposition~\ref{claim:nonneg} we already know that this resistance is non-negative, but more can be said for the cost node. Refer to Fig.~\ref{fig:reducedNetw} showing a subnetwork connecting two arbitrary nodes. When one of the nodes is the cost node, it does not have one of the negative resistances, therefore, the the total resistance, $R_{\text{total}}$, which can be seen from this node is at least all the cost resistances in parallel
\begin{align}
R_{\text{total}} \geq \frac{1}{\sum_{i=1}^n {\Gcs}_i}.
\label{eq:Rtot}
\end{align}
From~\eqref{eq:fcUIcost} follows that
\begin{align}
& \Gcs^T \Delta V = \left (\sum_{i=1}^n {\Gcs}_i  \right )\Delta  \Ucs + \Delta  \Ics.
\label{eq:dCost}
\end{align}
Using the total equivalent resistance we know that
\begin{align}
& \Delta  \Ics  = - \frac{\Delta  \Ucs }{R_{\text{total}}} \label{eq:dIcs}.
\end{align}
Combination of~\eqref{eq:dCost},~\eqref{eq:dIcs} and~\eqref{eq:Rtot} yields
\begin{align}
& \frac{\Gcs^T \Delta V}{\Delta  \Ucs} = \sum_{i=1}^n {\Gcs}_i  - \frac{1}{R_{\text{total}}} \geq 0.
\label{eq:DcostdUcs}
\end{align}
The equation~\eqref{eq:DcostdUcs} states that the change in cost value must have the same sign as the change in $\Delta  \Ucs$. Therefore, when $\Ucs$ is decreased the cost must decrease or stay the same. However, the cost cannot decrease, since it is already optimal. Therefore the cost must remain constant, and the circuit holds solution to the problem~\eqref{eq:genopt} for any $\Ucs \leq \Ucs^{crit}$.
%
This result completes the proof of Theorem~\ref{thrm:equiv}.



\section{Example applications and experimental results}
\label{sec:examp}
This section presents three examples where the approach proposed in this paper has been successfully applied.
In the first example an LP is solved by the proposed  electrical circuit simulated by using the SPICE~\cite{SPICE} simulator.
In the second example an analog LP is used to control a linear system by using Model Predictive Control.
In the third example  an experiment is conducted by realizing the circuit for a small LP with standard electronic components.

\subsection{Linear Programming}
We demonstrate capability of the method by solving an LP problem. The problem is a randomly generated and  it has 120 variables, 70 equality constraints and 190 inequality constraints.
In order to simulate parasitic effects of real circuit inductance values of $100 nH$ are assumed for the wires, that roughly corresponds to inductance of 10~cm long wire.

\begin{figure}[tb]
\centering
   \includegraphics[width=0.35\textwidth, trim=0.1cm 0 0.85cm 1.4cm, clip]{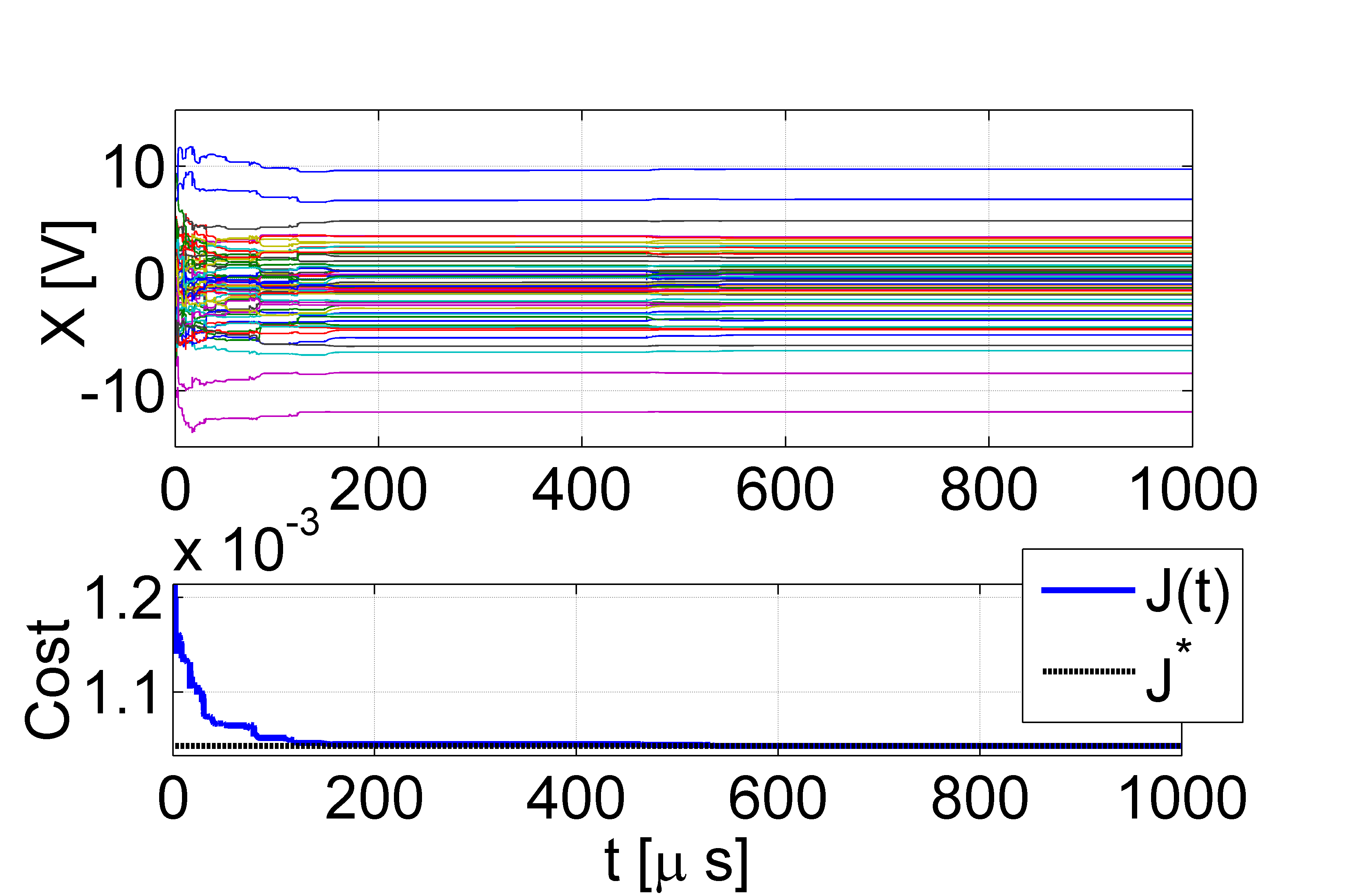}\\
  \caption{Example of LP solution. The upper plot shows solution variables in time. The lower plot shows the cost function value.}\label{fig:LPex}
\end{figure}

The convergence of the electric circuit is shown in Fig.~\ref{fig:LPex}.
The time scale in this example is determined by the selected value of parasitic inductance. The circuit transient can be partitioned to two phases.
During the first $200 \mu s$ rapid convergence  to a solution close to the optimal one can be observed. Afterwards, at about $500 \mu s$ the circuit converges to the true optimum value. Typical accuracy achieved in analog electronics is in the order of $0.5\%$ of the dynamic range.
The longer convergence time is not of practical interest, because the difference between the immediate cost value and the true optimal one is less than the accuracy that is expected from analog devices.

\subsection{MPC example}
This example demonstrates the implementation of a model predictive controller with an LP analog circuit.
For this example we work with the dynamical system
$\frac{dx}{dt} = -x + u$,
where $x$ is the system state and $u$ is the input.
We want $x$ to follow a given reference trajectory, while satisfying input constraints.
The finite time optimal control problem at time $t$ is formulated as
\begin{subequations}
\label{eq:mpcex}
\begin{align}
&\min_{u_{0}\ldots u_{n-1}} \sum_{i=1}^N | x(i)- x_{ref}(i) | \\
& x_{i+1} = x_i +  (u_i-x_i)\delta, \ i=0,\ldots,N \\
& -1.5 \leq u_i \leq 1.5, \ i=0,\ldots,N \\
& x_0 = x(t)
\end{align}
\end{subequations}
where $N$ is the prediction horizon, $x_{ref}(i)$ is the reference trajectory at step $i$, $\delta$ is sampling time  and $x(t)$ is the initial state at time $t$. Only the first input, $u_0$, is applied at each time step $t$.

\begin{figure}[tb]
\centering
   \includegraphics[width=0.4\textwidth, trim=0.5cm 0 1.5cm 0cm, clip]{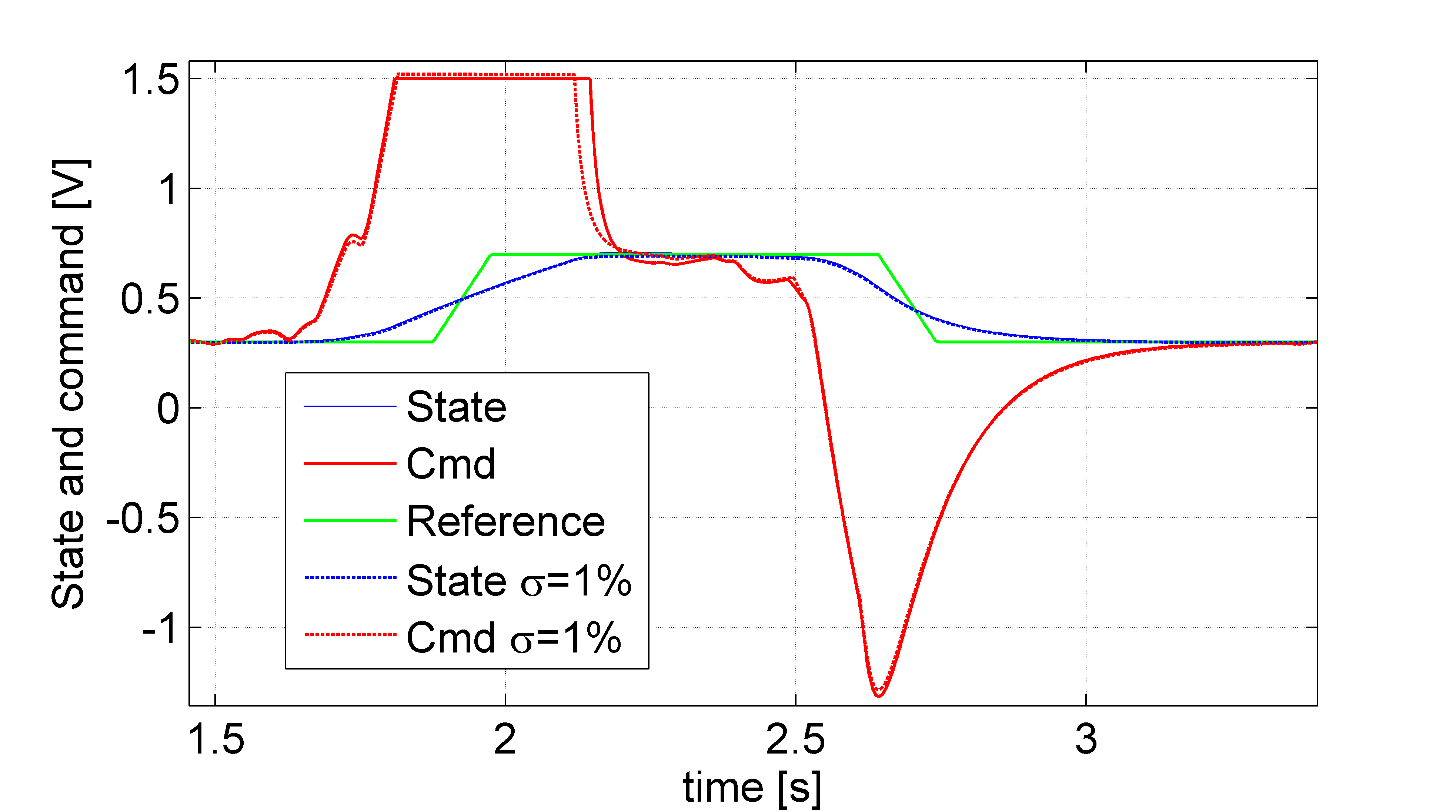}\\
  \caption{Example of MPC implementaion. Solid lines represent nominal controller, dashed lines represent controller implemented with random $1\%$ error of analog devices.}\label{fig:MPCex}
\end{figure}

With $N=16$, the LP in~\eqref{eq:mpcex} has 96 variables, 63 equality constraints and 49 inequality constraints.
An electric circuit that implements system dynamics together with the circuit that implements the MPC controller were constructed and simulated using SPICE. The voltage value representing the system state was measured and enforced on the $x_0$ node of the LP. The optimal input value $u_0$ was injected as input to the simulated system dynamics.
Fig.~\ref{fig:MPCex} shows the closed loop simulations results. Notice the predictive behavior of the closed loop control input and the satisfaction of the system constraints.

%

In order to demonstrate system performance for imperfect analog devices, another simulation result with $1\%$ random Gaussian error in values of resistors is presented on the same Fig.~\ref{fig:MPCex}. There is no significant change in system behavior.

\subsection{Hardware implementation example}

We implemented a small LP using standard electronics components. The same problem was realized   by Hopfield~\cite{Hopfield86} and Chua~\cite{Chua88}. The LP is defined as follows
\begin{eqnarray}
\min_{x_1,x_2} c^T [x_1 \ x_2]^T & \notag \\
s.t. \  \frac{5}{12}x_1   - x_2 \leq \frac{35}{12} , \
& \frac{5}{2}x_1   + x_2  \leq  \frac{35}{2} \notag \\
 -x_1  \leq 5 , \
& x_2   \leq 5
\end{eqnarray}
where $c$ is a cost vector, that is varied to get different solution
points. The circuit was realized using resistors of $1\%$ accuracy, operational
amplifiers (OP27) for the negative resistance and comparator (LM311) together
with the switch (DG201) to implement functionality of an ideal diode .

Various values for the cost function $c$ and test results are summarized in
Table~\ref{tbl:hardresults}. Table~\ref{tbl:hardresults} shows that the experimental results are
accurate up to 0.5\%. The circuit reaches an equilibrium  $6~\mu s$ after
the cost voltage was applied. The convergence time is governed by a slew rate
of the OP27 that is limited to 2.8 $V/\mu s$.
\begin{table}[tb]
\centering
\caption{Experimental and theoretical results (in parenthesis) for LP solution.   }
\label{tbl:hardresults}
\begin{tabular}{ccc}
\hline \hline
cost direction & x1 (exact) & x2 (exact) \\
\hline
1 1 & 4.996 (5.0)& 4.99 (5.0)\\
-1 1 & 7.002 (7.0)& 5.005 (5.0)  \\
-1 -1 & -7.012 (-7.0)& -4.98 (-5.0) \\
1 0 & 6.976 (7.0) & 0.005 (0.0) \\
\hline
\end{tabular}
\end{table}

\section{Conclusion}
\label{sec:concl}

In this paper we presented an approach to design an electric analog circuit that
is able to solve a feasible Linear Program.
The method is used to implement and solve MPC based on linear programming.
We present simulative and the experimental results that demonstrate the effectiveness of the proposed method.

The reported LP solution speed of $6~\mu s$  is faster than any result that was previously reported in the literature, and may be
significantly decreased further by selecting faster components or
implementing the design using faster technology, such as custom VLSI design or
FPAA device.

The circuit analysis is at steady state. The theory of Linear Complimentary system~\cite{Heemels98} can be used to study the dynamic circuit behavior. This is a subject of ongoing research.
Future research directions include solution of larger problems, possible expansion the method to solution of quadratic programming
(QP) and solutions to the optimal circuit design.

\section{Acknowledgments}
The authors would like to thank prof. Ilan Adler for valuable discussions that provided helpful inputs to this work, such as the primal-dual LP circuit. Also we gratefully acknowledge the financial support of Helen Betz Foundation for this research.

\bibliographystyle{ieeetr}
\bibliography{bibliography}

\end {document}